\documentstyle[11pt]{amsart}

\newtheorem{prop}{Proposition}
\newtheorem{theo}{Theorem}
\newtheorem{Lemma}{Lemma}
\newtheorem{cor}{Corollary}

\newcommand{\na}{\nabla}
\newcommand{\om}{\omega}
\newcommand{\Om}{\Omega}
\newcommand{\la}{\lambda}
\newcommand{\La}{\Lambda}

\newcommand{\ka}{K{\"a}hler }

\newtheorem{conjecture*}{Conjecture}
\newtheorem{theorem*}{Theorem}
\newtheorem{question*}{Question}

\newcommand{\leftr}{[\hbox{\hspace{-0.15em}}[}
\newcommand{\rightr}{]\hbox{\hspace{-0.15em}}]}

  \newcommand{\lllra}{\!\!\joinrel{\hbox to
50pt{\rightarrowfill}}} \newcommand{\llllra}{\!\!\joinrel{\hbox to
60pt{\rightarrowfill}}} \newcommand{\lllllra}{\!\!\joinrel{\hbox to
70pt{\rightarrowfill}}} \newcommand{\llllllra}{\!\!\joinrel{\hbox to
75pt{\rightarrowfill}}}

\title[The curvature and the integrability of almost-K{\"a}hler
manifolds]{The curvature and the integrability of almost-K{\"a}hler
manifolds: a survey}
\author{VESTISLAV APOSTOLOV AND TEDI DR\u{A}GHICI} \thanks{The first author
was supported in part by a PAFARC-UQAM grant, NSERC grant OGP0023879 and by
FCAR grant NC-7264. He is also member of EDGE, Research Training Network
HPRN-CT-2000-00101, supported by the European Human Potential Programme.
Both authors were also supported in part by NSF grant INT-9903302}
\address{Vestislav Apostolov \\ D\'epartement de math\'ematiques\\ UQAM \\
succursale Centre-ville c.p. 8888 \\ Montr\'eal \\ H3C 3P8, Canada}
\email{apostolo@@math.uqam.ca}

\address{Tedi Dr\u{a}ghici \\ Department of Mathematics \\ Florida
International University \\ Miami FL 33199 \\ USA}
\email{draghici@@fiu.edu}

\begin{document}

\maketitle

\section{Introduction}

Situated at the intersection of Riemannian, symplectic and complex
geometry, K\"ahler geometry is a subject that has been intensively studied
in the past half-century. The overlap of just Riemannian and symplectic
geometry is thewider area of almost-\ka geometry. Although problems
specific to almost-\ka geometry have been considered in the past (see e.g.
\cite{Go,gray2,olszak}), it is fair to say that our knowledge beyond the
\ka context is still small at this point. In recent years there is a
renewed interest in almost-\ka geometry, mainly motivated by the rapid
advances in the study of symplectic manifolds.

\vspace{0.2cm} An {\it almost-\ka structure} on a manifold $M$ of (real)
dimension $n = 2m$ is a blend of three components: a symplectic form
$\Omega$, an almost-complex structure $J$ and a Riemannian metric $g$,
which satisfy the compatibility relation \begin{equation} \label{compat}
\Om(\cdot, \cdot) = g(J\cdot, \cdot).
\end{equation}
Any one of these three components is completely determined by the remaining
two via (\ref{compat}). If, additionally, the almost-complex structure $J$
is integrable, i.e., is induced by a complex coordinate atlas on $M$, then
the triple $(g, J, \Om)$ is a {\it K\"ahler} structure on $M$. In real
dimension 2, the notions of almost-\ka and \ka structure coincide, but this
is no longer true in higher dimensions. Throughout the paper {\it strictly}
almost-K\"ahler will mean that the corresponding almost-complex structure
is {\it not} integrable, equivalently, that the almost-\ka structure is not
K\"ahler.

\vspace{0.2cm}

Depending on whether the Riemannian metric, or the symplectic form are
initially given, there are slightly different perspectives towards studying
almost-\ka structures. Given a Riemannian metric $g$, it is a difficult
question to decide if $g$ admits any compatible almost-\ka structure at
all; in other words, we ask for the existence of an almost-complex
structure $J$ such that $\Om$ determined by (\ref{compat}) is a closed
2-form. It is well known
that the holonomy group determines if the given metric admits a compatible
\ka structure. Comparatively, the almost-\ka problem is considerably more
intricate, even locally. For example, until recently it was still an open
question whether or not metrics of constant negative sectional curvature
admit compatible almost-\ka structures. Now the answer is known to be
negative as a consequence of the results in
\cite{Arm2,oguro-sekigawa,olszak}.

\vspace{0.2cm}

Starting with a symplectic manifold $(M, \Om)$, there are many pairs
$(g,J)$ which satisfy the compatibility relation (\ref{compat}) and we
shall refer to any such pair as an $\Om$-compatible almost-\ka structure on
$M$. The space of all $\Om$-compatible almost \ka structures ---
equivalently, the space of all $\Om$-compatible metrics (or almost-complex
structures) --- will be denoted by ${\bf AK}(M,\Om)$ and is well known to
be an infinite dimensional, contractible Fr{\'e}chet space. The choice of
an {\it arbitrary} (or generic) metric in this space serves as an useful
tool in studying the symplectic geometry and topology of the manifold. This
is manifested in the theory of pseudo-holomorphic curves initiated by
Gromov \cite{gromov}, and more recently, in Taubes' characterization of
Seiberg-Witten invariants of symplectic 4-manifolds \cite{taubes1,taubes2}.
Rather than dealing with generic $\Om$-compatible almost-\ka metrics, the
present survey is centered around the problem of identifying and studying
``{\it distinguished}'' Riemannian metrics in ${\bf AK}(M,\Om)$.

\vspace{0.2cm}

One particularly nice situation would be to have an $\Om$-compatible \ka
metric $g$. In this case, the Levi-Civita connection $\nabla$ preserves the
other two structures, i.e., $\na J =0$ and $\na \Om =0$. This leads to
useful symmetry properties of basic geometric operators (like Laplacian,
curvature, etc.) with respect to the complex structure, which further
determine geometric and topological consequences on the manifold. However,
it is now well known that in dimension higher or equal to four, most
compact symplectic manifolds do not admit \ka metrics. Thus one should
require less of a compatible metric.

\vspace{0.2cm} A good candidate for privileged Riemannian metric on a given
manifold is an {\it Einstein}metric, i.e., a Riemannian metric $g$ for which
\begin{equation}\label{einstein}
{\rm Ric} = \lambda g,
\end{equation}
where ${\rm Ric}$ is the Ricci tensor and $\la$ is a constant equal, up to
a factor $\frac{1}{{\rm dim}_{\mathbb R}(M)}$, to the scalar curvature of
$g$. Although a remarkable amount of research in Riemannian geometry has
been done to study Einstein metrics (see e.g. \cite{besse,newbesse}), it
seems that the almost-K\"ahler aspects of the problem remain in its
infancy. A long standing, still open conjecture of Goldberg \cite{Go}
affirms that there are no {\it Einstein}, strictly almost-\ka metrics on a
compact symplectic manifold. Indirectly, the Goldberg conjecture predicts
that compatible Einstein metrics are very rare on compact symplectic
manifolds. The conjecture is still far from being solved, but there are
cases when it has been be confirmed: Sekigawa \cite{sekigawa} proved that
the conjecture is true if the scalar curvature is non-negative (see Theorem
\ref{thm:sekigawa} below) and there are further positive partial results in
dimension four under other additional curvature assumptions
\cite{AA,Arm1,Arm2,oguro-sekigawa}. Moreover, a number of subtle
topological restrictions to the existence of Einstein metrics on compact
4-manifolds are now known \cite{kotschick,lebrun1,lebrun2,lebrun3}, and
these can be thought as further support for the conjecture.

\vspace{0.2cm}

Although the Goldberg conjecture is global in nature, it is the merit of
John Armstrong \cite{Arm0} to point out that the local aspects of the
problem will most likely determine the global ones as well. His idea was to
apply Cartan-\ka theory to prove local existence of strictly almost-\ka
Einstein structures. Because of mounting algebraic difficulties he did not
achieve this goal, but he proved several other interesting results and it
became clear that Cartan-\ka theory is a very useful tool for various
problems in almost-\ka geometry. For example, it turns out that the
Cartan-\ka theory is not uniquely relevant to Einstein almost-\ka metrics,
but is also suited for the study of a more general class of almost-\ka
metrics arising from a variational approach \cite{blair,BI}. This approach
is reviewed in section 2. After we introduce some notations and basic facts
of almost-\ka geometry in section 3, section 4 casts some of the known
integrability results under the unifying framework of Cartan-\ka theory. In
section 5 we then present a number of local examples of strictly almost-\ka
manifolds with interesting curvature properties, in particular some recent
local constructions of Einstein strictly almost-\ka metrics. The material
in section 6 illustrates how additional conditions on the curvature often
lead to completely solvable equations. We present some recent local
classification results in dimension four, which can be used to obtain
further integrability results. The last section is inspired by the recent
works of Donaldson \cite{donaldson1} and C. LeBrun \cite{lebrun3}, which
naturally lead to the study of Einstein-like conditions with respect to the
canonical Hermitian connection. We present a brief discussion on the
relevant curvature conditions and highlight some directions for further
research.

\section{A variational approach to almost-K{\"a}hler manifolds}

Recall that \cite{besse} on a compact manifold $M$, the Einstein condition
is the Euler-Lagrange equation of the {\it Hilbert functional} ${\bf S}$,
the integral of the scalar curvature, acting on the space of all Riemannian
metrics on $M$ of a given volume.

A ``symplectic'' setting of this variational problem was proposed by Blair
and Ianus \cite{BI}: restricting the Hilbert functional to the space ${\bf
AK}(M,\Om)$ of compatible metrics of a given compact symplectic manifold
$(M,\Omega)$, the critical points are the almost-\ka metrics $(g,J)$ whose
Ricci tensor ${\rm Ric}$ is $J$-invariant, i.e., satisfies
\begin{equation}\label{eq:1}
{\rm Ric}(J\cdot, J\cdot) = {\rm Ric}(\cdot,\cdot). \end{equation} The
Euler-Lagrange equation (\ref{eq:1}) is thus a weakening of both the
Einstein and the \ka conditions. Furthermore, Blair \cite{blair} observed
that for any almost-\ka metric $(g,J)$ the following relation holds:
\begin{equation}\label{eq:blair}
\frac{1}{4} \int_M |\nabla J|^2 dv + {\bf S}(g) = \frac{4\pi}{(m-1)!} (c_1
\cdot [\Omega]^{\wedge (m-1)})(M), \end{equation} where $\nabla$ is the
Levi-Civita connection of $g$, $|\cdot |$ is the point-wise norm induced by
$g$, and $c_1$ and $dv = \frac{1}{m!}\Omega^{\wedge m}$ are respectively
the first Chern class and the volume form of $(M,\Om)$. It follows that
${\bf S}$ is directly related to the {\it Energy functional} which acts on
the space of $\Om$-compatible almost-complex structures by
$${\bf E} (J) = \int_M |\na J|^2 dv .$$
>From this point of view, almost-\ka metrics satisfying
(\ref{eq:1}) have been recently studied in \cite{Le-Wang}, and have been
called there {\it harmonic} almost-\ka metrics. However, in this paper we
adopt the definition:

\vspace{0.2cm} \noindent {\it Definition 1.} An almost-\ka metric $(g,J)$
is called {\it critical} if it satisfies (\ref{eq:1}).

\vspace{0.2cm} \noindent It follows from (\ref{eq:blair}) that the
functional ${\bf E}$ (resp.~${\bf S}$) is bounded from below (resp. from
above), the \ka metrics being minima of ${\bf E}$ (resp. maxima of ${\bf
S}$). However, a direct variational approach to finding extrema seems not
to be easily
applicable since it may happen that the infimum of ${\bf E}$ be zero,
although $M$ does not carry \ka structures at all.

\vspace{0.2cm} \noindent {\it Example 1.} \cite{Le-Wang} Let $M = S^1\times
(Nil^3/\Gamma)$, where $Nil^3$ is the three-dimensional Heisenberg group
and $\Gamma$ is a co-compact lattice of $Nil^3$; as first observed by
Thurston \cite{thurston}, $M$ is a smooth four-dimensional manifold which
admits a symplectic structure but does not admit \ka metrics at all.
Specifically, $M$ carries an invariant complex structure $J$ with trivial
canonical bundle, and thus comes equipped with a holomorphic symplectic
structure. To see this we write
$$Nil^3 = \{ A \in GL(3,{\mathbb R}) | \; A = \left( \begin{array}{ccc} 1 &
x & z \\
0 & 1 & y \\
0 & 0 & 1 \end{array} \right), \; x, y, z \in {\mathbb R} \}.$$ Then the
1-forms $dt, dx, dy, dz-xdy$ are left-invariant on $S^1\times Nil^3$; we
define a left-invariant complex structure $J$ on this manifold by $$J(dx)
:= dy; \ J (dz-xdy) : = dt, $$
and a holomorphic symplectic form $\om$ by $$\om = \Big(dx\wedge(dz - xdy)
- dy\wedge dt \Big) +i\Big(dx\wedge dt + dy\wedge(dz-xdy) \Big).$$ Being
left-invariant, both $J$ and $\om$ descend to $M$ to define a so-called
Kodaira complex surface $(M,J,\om)$, cf. e.g. \cite{BPV}. Note that the
first Betti number of $M$ is equal to $3$, showing that no \ka metric
exists on $M$. On the other hand
$$ \Om = {\rm Im}(\om) = dx\wedge dt + dy\wedge(dz-xdy)$$ is a
left-invariant symplectic form. Following \cite{Le-Wang}, we consider a
2-parameter family ${\bf AK}(a,b), a>0, b>0,$ of $\Om$-compatible
(left-invariant) almost-\ka structures $(g_{a,b}, I_{a,b}, \Om)$, defined
by $$g_{a,b} = \frac{1}{a^2}dx^{\otimes 2} + \frac{1}{b^2}dy^{\otimes 2} +
a^2 dt^{\otimes 2} + b^2(dz-xdy)^{\otimes 2};$$ $$I_{a,b}(dx) := a^2dt; \
I_{a,b}(dy) := b^2(dz-xdy).$$ One calculates directly $${\rm Ric}(g_{a,b}) =
-\frac{a^2b^4}{2}\Big(\frac{1}{a^2}dx^{\otimes 2} +
\frac{1}{b^2}dy^{\otimes 2} - b^2(dz-xdy)^{\otimes 2}\Big),$$ so that the
scalar curvature is
$$s_{a,b}= -\frac{a^2b^4}{2}.$$
It follows that on the space ${\bf AK}_{a,b}$ we have ${\rm sup} ({\bf H})
= {\rm inf}({\bf E}) =0$, while $M$ does not carry \ka metrics at all.
$\Diamond$

\vspace{0.2cm}

In general, for a given compact symplectic manifold $(M^{2m},\Om)$, it is
tempting to compute the symplectic invariant $$ {\bf E} (\Omega) = \inf_{g}
\int_M |\na J|^2 d\mu, $$ where the infimum is taken over all almost-\ka
structures $(g,J)$ compatible with $\Omega$. A natural approach to this
problem is to study the (global) behavior of the solutions of the
corresponding (negative) {\it gradient flow} equation $$ \frac{d}{dt} g_t =
- ({\rm Ric}^t)'', $$ where $(g_t,J_t)$ is a path of almost-\ka structures
in ${\bf AK}(M,\Om)$, ${\rm Ric}^t$ is the Ricci tensor of $g_t$ and $({\rm
Ric}^t)''$ is the gradient of ${\bf E}$, i.e., $$ ({\rm
Ric}^t)''(\cdot,\cdot) = \frac{1}{2}\Big({\rm Ric}^t(\cdot,\cdot) - {\rm
Ric}^t(J_t \cdot, J_t \cdot) \Big).$$ This approach has been recently
explored in \cite{Le-Wang}. However, as Example 1 suggests, technical
difficulties appear when studying global properties of the flow $g_t$. One
reason for this is perhaps hidden in the fact that even locally the
equation (\ref{eq:1}) is difficult to be solved, and our next goal is to
highlight some of these ``local obstructions''.

\section{The curvature and the Nijenhuis tensor: First obstructions}
\subsection{Type-decompositions of forms and vectors.}

Let $(M,g,J,\Om)$ be an almost-\ka manifold of real dimension $n=2m$. We
denote by: $TM$ the (real) tangent bundle of $M$; $T^*M$ the (real)
cotangent bundle; $\La^r M, r=1,...,n$ the bundle of real $r$-forms;
$S^{\ell}M$ the bundle of symmetric $\ell$-tensors; $(\cdot,\cdot )$ the
inner product induced by $g$ on these bundles (or on their tensor
products).

Using the metric, we shall implicitly identify vectors and covectors and,
accordingly, a 2-form $\phi$ with the corresponding skew-symmetric
endomorphism of the tangent bundle $TM$, by putting: $( \phi(X),Y ) =
\phi(X,Y)$ for any vector fields $X,Y$.

The almost-complex structure $J$ gives rise to a type decomposition of
complex vectors and forms. By convention, $J$ acts on the cotangent bundle
$T^*M$ by $(J\alpha)_X = -\alpha_{JX}$, so that $J$ commutes with the
Riemannian duality between $TM$ and $T^*M$. We shall use the standard
decomposition of the complexified cotangent bundle $$T^*M\otimes {\Bbb C}=
\La^{1,0}M \oplus \La^{0,1}M,$$ given by the $(\pm i)$-eigenspaces of $J$,
the type decomposition of complex $2$-forms
\begin{equation*}\label{Lambda2c}
\La^2M\otimes {\mathbb C} = \La^{1,1}M \oplus \La^{2,0}M \oplus \La^{0,2}M,
\end{equation*}
and the type decomposition of symmetric (complex) bilinear forms $$S^2M
\otimes {\mathbb C} = S^{1,1}M \oplus S^{2,0}M \oplus S^{0,2}M.$$ Besides
the type decomposition of complex forms, we shall also consider the
splitting of real 2-forms into $J$-invariant and $J$-anti-invariant ones;
for any real section of $\La^2M$ (resp. of $S^2 M$), we shall use the
super-script $'$ to denote the projection to the real sub-bundle
$\La^{1,1}_{\Bbb R} M$ (resp. $S^{1,1}_{\Bbb R} M$) of $J$-invariant
2-forms (resp. symmetric 2-tensors), while the super-script $''$ stands for
the projection to the bundle $\leftr \La^{0,2}M \rightr $ (resp. $\leftr
S^{0,2}M \rightr$) of $J$-anti-invariant ones; here and henceforth $\leftr
\ \cdot \ \rightr$ denotes the real vector bundle underlying a given
complex bundle. Thus, for any section $\psi$ of $\La^2 M$ (resp. of $S^2
M$) we have the orthogonal splitting $\psi = \psi' + \psi'' ,$ where
$$\psi'(\cdot , \cdot) = \frac{1}{2}(\psi(\cdot , \cdot) + \psi(J\cdot,
J\cdot)) \; \mbox{ and } \; \psi''(\cdot , \cdot) = \frac{1}{2}(\psi(\cdot
, \cdot) - \psi(J\cdot, J\cdot)) .$$ The real bundle $\leftr \La^{0,2}M
\rightr$ (resp. $\leftr S^{0,2}M\rightr $) inherits a canonical complex
structure, still denoted by ${J}$, which is given by
$$({J}\psi)(X,Y) := -\psi(JX,Y),
\ \forall \psi \in \leftr \La^{0,2}M \rightr, $$ so that $(\leftr
\La^{0,2}M \rightr, J)$ becomes isomorphic to the complex bundle
$\La^{0,2}M$. We adopt a similar definition for the action of ${J}$ on
$\leftr S^{0,2}M\rightr $. Notice that, using the metric $g$, $\leftr
S^{0,2}M\rightr$ can be also viewed as the bundle of symmetric,
$J$-anti-commuting endomorphisms of $TM$.

We finally define the $U(m)$-decomposition (with respect to $J$) of real
2-forms
\begin{equation}\label{Lambda2r}
\La^2M = {\mathbb R}\cdot \Om \oplus \La^{1,1}_0 M \oplus \leftr \La^{0,2}M
\rightr,
\end{equation}
where $\La^{1,1}_0 M$ is the sub-bundle of the {\it primitive} (1,1)-forms,
i.e., the $J$-invariant 2-forms which are point-wise orthogonal to $\Om$.

\subsection{Decompositions of algebraic curvature tensors} Following
\cite{besse}, we denote by
$R$ the curvature with respect to the Levi-Civita connection $\na$ of $g$.
Classically, $R$ is considered as a real $(1,3)$-tensor on $M$, but, using
the metric, we shall rather think of $R$ as a $(0,4)$-tensor, or as a
section of the bundle $S^2(\La^2M)$ of symmetric endomorphisms of $\La^2M$,
depending on the context. We adopt the convention $R_{XYZV} = (\na_{[X,Y]}Z
- [\na_X, \na_Y]Z,V);$ the action of $R$ on $\La^2 M$ is then defined by
$$R(\alpha\wedge \beta)_{X,Y} := R_{\alpha^{\sharp}, \beta^{\sharp},X,Y},$$
where the super-script $\sharp$ denotes the corresponding (Riemannian) dual
vector fields. Thus, the curvature $R$ can be viewed as a section of the
symmetric tensor product $\La^2M\odot \La^2M$, and satisfies the {\it
algebraic Bianchi identity} $$ R_{XYZ} + R_{ZXY} + R_{YZX}=0.$$
It follows that $R$ is actually
a section of the sub-bundle
${\cal R}(M) = {\rm ker} ({\frak e})$, where ${\frak e} : \La^2M\odot
\La^2M \to \La^4 M$ is given by ${\frak e}(\phi \odot \phi) =
\phi\wedge\phi$. The bundle ${\cal R}(M)$ is called the bundle of {\it
algebraic curvature tensors}.

Recall that the {\it Ricci} contraction, ${\rm Ric}$, associates to $R$ the
symmetric bilinear form $${\rm Ric}(X,Y) = \sum_{i=1}^{n}
R_{X,e_i,Y,e_i}.$$ It can be thought as a linear map ${\frak r} : {\cal
R}(M) \to S^2M$. For $n\ge 3$ this map is surjective and its adjoint,
${\frak r}^*$, is injective. We thus obtain the orthogonal decomposition
$${\cal R}(M) = {\frak r}^*(S^2M) \oplus {\cal W}(M),$$ where the
sub-bundle ${\cal W}(M)$ is the kernel of ${\frak r}$ and is called the
bundle of {\it algebraic Weyl tensors}. Accordingly, the curvature $R$
splits as $R ={\widehat {\rm Ric}} + W$, where the {\it Weyl curvature},
$W$, is the component of $R$ in ${\cal W}(M)$ while ${\widehat {\rm Ric}}$
is identified with the Ricci tensor by the property ${\frak r}({\widehat
{\rm Ric}}) = {\rm Ric}$. Specifically, starting from the splitting of
${\rm Ric}$ into its trace $$s= \sum_{i=1}^{n} {\rm Ric}(e_i,e_i),$$ the
scalar curvature of $(M,g)$, and its trace-free part ${\rm Ric}_0 = {\rm
Ric} - \frac{s}{n}g$, we have that $${\widehat {\rm Ric}} = h\wedge g, $$
where $h= \frac{s}{2n(n-1)}g + \frac{1}{(n-2)}{\rm Ric}_0$ is the {\it
normalized Ricci tensor} of $(M,g)$ and $\wedge$ denotes the
Kulkarni-Nomizu product of bilinear forms. We finally obtain the splitting
of $R$ into three orthogonal pieces: \begin{equation}\label{o(n)} R =
\frac{s}{n(n-1)}{\rm Id}|_{\La^2M} + \frac{1}{(n-2)}({\rm Ric}_0\wedge g) +
W. \end{equation}
Each of the three terms in the right-hand side of (\ref{o(n)}) is an
element of ${\cal R}(M)$; the corresponding sub-bundle of ${\cal R}(M)$ is
attached to an irreducible representation of the orthogonal group $O(n)$
acting on ${\cal R}({\mathbb R}^n)$, namely the trivial representation, the
Cartan product ${\mathbb R}^n \odot {\mathbb R}^n$ and the product ${\frak
o}(n)\odot {\frak o}(n)$ where the Lie algebra ${\frak o}(n)$ is identified
with $\La^2 {\mathbb R}^n$. Note that the {\it Einstein equations} for a
Riemannian metric $g$ correspond to the vanishing of the component
$\frac{1}{(n-2)}({\rm Ric}_0\wedge g)$ of $R$; equivalently, the Ricci
tensor satisfies ${\rm Ric} = \frac{s}{n}g$, in which case the scalar
curvature $s$ is necessarily constant. Recall also that a manifold with
zero Weyl tensor is conformally flat.

In dimension four, one can further refine the splitting (\ref{o(n)}) by
fixing an orientation on $(M,g)$ and considering the {\it Hodge}
$*$-operator which acts on the bundle $\La^2M$ as an involution, thus
giving rise to an orthogonal splitting of $\La^2M$ into the sum of the
$\pm$-eigenspaces, the bundles of self-dual and anti-self-dual 2-forms.
Consequently,
the Weyl tensor also decomposes as $W^+ + W^-$, where $W^+ = \frac{1}{2}(W
+ W\circ *)$ is the part acting on self-dual forms, while $W^- =
\frac{1}{2}(W - W \circ *)$ is the restriction of $W$ to the bundle of
anti-self-dual forms. An oriented Riemannian 4-manifold $(M,g)$ is called {\it self-dual} (resp. {\it anti-self-dual}) if its Weyl tensor is self-dual (re
sp. anti-self-dual), i.e., if $W^- \equiv 0$ (resp. $W^+ \equiv 0$).

We now assume that $n=2m$ is even and that $(M^{2m},g)$ is endowed with a
$g$-orthogonal almost-complex structure $J$, i.e., $(g,J)$ is an
almost-Hermitian structure on $M$. Representation theory of the unitary
group $U(m)$ implies the splitting of the bundle ${\cal R}({\mathbb C}^m)$
into irreducible factors, each of them defining a {\it curvature component}
of $R$ with respect to $J$, cf. \cite{TV}. Thus, for $n \ge 8$, $R$ has
$10$ different curvature components; for $n=6$, we have $9$ curvature
components, while for $n=4$ one gets $7$ different curvature components.
Note that the most natural special types of almost-Hermitian manifolds
correspond to the vanishing of some of these components.

The
following $U(m)$-invariant sub-bundles of ${\cal R}(M)$ arise naturally
\cite{gray2,TV}:
\begin{enumerate}
\item[$\bullet$] ${\cal K}(M) = \{ R \in {\cal R}(M) : R_{JXJY} = R_{XY}
\}$ is the sub-bundle of elements acting trivially on $\leftr \La^{2,0}M
\rightr$. The curvature of any \ka metric belongs to ${\cal K}(M)$, so
${\cal K}(M)$ is usually referred to as the bundle of {\it k\"ahlerian
algebraic curvature tensors}.

\item[$\bullet$] ${\cal R}_2(M) = \{ R \in {\cal R}(M) : R_{XYZW} -
R_{JXJYZW} = R_{JXYJZW} + R_{JXYZJW} \}$ is the sub-bundle of elements
preserving the type decomposition (\ref{Lambda2c}) of complex 2-forms;
equivalently, these are the elements of ${\cal R}(M)$ which preserve the
splitting (\ref{Lambda2r}) and commute with the complex structure of
$\leftr \La^{0,2}M \rightr$.

\item[$\bullet$] ${\cal R}_3(M) = \{ R \in {\cal R}(M) :
R_{JXJYJZJT}=R_{XYZT} \}$ is the
sub-bundle of elements of ${\cal R}(M)$ preserving the type decomposition
(\ref{Lambda2r}) of real 2-forms. \end{enumerate}

Note that ${\cal K}(M) \subset {\cal R}_2(M) \subset {\cal R}_3(M)$, and
for any almost-Hermitian manifold $(M,g,J)$ whose curvature tensor belongs
to
${\cal R}_3(M)$, the Ricci tensor satisfies (\ref{eq:1}).

Each of the
above sub-bundles can be further decomposed into irreducible factors. For
instance, the bundle
${\cal K}(M)$ splits into three pieces, associated respectively to the
trivial representation of
$U(m)$, to the Cartan product ${\mathbb C}^m\odot {\mathbb C}^m$, and to
${\frak u}(m)\odot {\frak u}(m)$. Corresponding to the trivial factor is
the curvature of a complex space form (so that this component can be
identified to the scalar curvature), the second factor is identified to the
traceless Ricci tensor, while the third component is the so-called {\it
Bochner curvature}; for more details see \cite{TV,FFS}.

We further denote by
$W''$ the component of the curvature operator defined by $$W''_{XYZT} =
\frac{1}{8}\Big(R_{XYZT} - R_{JXJYZT} - R_{XYJZJT} + R_{JXJYJZJT}$$ $$
\hspace{3cm} - R_{XJYZJT} - R_{JXYZJT} - R_{XJYJZT} + R_{JXYJZT}\Big).$$ It
follows that
$$W''_{Z_1Z_2Z_3Z_4} =
R_{Z_1Z_2Z_3Z_4}= W_{Z_1Z_2Z_3Z_4} \ \ \forall Z_i \in T^{1,0}M, $$ which
explains the notation.

Let us also introduce the
following Ricci-type curvature tensors of an almost-Hermitian manifold:
\begin{enumerate}
\item[$\bullet$] the invariant and the anti-invariant parts of the Ricci
tensor with respect to $J$, denoted by ${\rm Ric}'$ and ${\rm Ric}''$
respectively. We also put $\rho(\cdot, \cdot) = {\rm Ric}'(J\cdot, \cdot)$
to be the (1,1)-form corresponding to the $J$-invariant part of ${\rm
Ric}$, which will be called {\it Ricci form} of $(M,g,J)$.
\item[$\bullet $] the {\it twisted Ricci form} $\rho^* = R(\Om)$ which is
not, in general, $J$-invariant. We will consequently denote by $(\rho^*)'$
and $(\rho^*)''$ the corresponding projections onto the bundles
$\La^{1,1}_{\Bbb R} M $ and $\leftr \La^{0,2}M \rightr $, respectively.
Following \cite{TV}, the trace of $\rho^*$, $s^*=2(R(\Om), \Om)$, will be
called $*$-{\it scalar curvature}. \end{enumerate} For a K{\"a}hler
manifold $\rho = \rho^*$ is the usual Ricci form which is closed, of type
(1,1) and represents (up to a scaling factor $\frac{1}{2\pi}$) the de Rham
cohomology of the first Chern class of $(M,J)$. Neither of these properties
remains true for a generic almost-K{\"a}hler manifold.

\subsection{The Nijenhuis tensor}
The {\it Nijenhuis tensor} (or
{\it complex torsion}) of an almost-complex structure $J$ is defined by
$$N_{X,Y} = [JX,JY] - J[JX,Y] - J[X,JY] - [X,Y].$$ Thus, $N$ is a 2-form
with values in $TM$, which by Newlander-Nirenberg theorem \cite{Ne-Ni}
vanishes if and only if $J$ is integrable. Alternatively, the Nijenhuis
tensor can be viewed as a map from
$\La^{1,0}M$ to $\La^{0,2}M$ in the following manner: given a complex
(1,0)-form $\psi$, we define by $\partial \psi$ and ${\bar \partial} \psi$
the projectors of $d\psi$ to $\La^{2,0}M$ and $\La^{1,1}M$, respectively.
In general, $d\neq \partial + {\bar \partial}$, as $d\psi$ can also have a
component of type (0,2), which we denote by $N(\psi)$; writing $\psi =
\alpha +
iJ\alpha$ where $\alpha$ is a real 1-form, one calculates $$N(\psi)_{X,Y}=
\frac{1}{2}\alpha(N_{X,Y})=\frac{1}{4} \psi(N_{X,Y}), \ \ \forall X,Y \in
T^{0,1}M.$$

For an almost-\ka manifold the
vanishing of $N$ is equivalent to $\Om$ being parallel with respect to the
Levi-Civita connection $\na$ of $g$; specifically, the following identity
holds (cf. e.g. \cite{KN}):
\begin{equation*} \label{naOm-N}
(\na_X\Om)(\cdot,\cdot) = \frac{1}{2} ( JX, N(\cdot,\cdot) ). \end{equation*}

Using the above relation, we shall often think of $N$ as a $T^*M$-valued
2-form, or as a $\Lambda^2 M$-valued 1-form, by tacitly identifying $N$
with $\na \Om$. Further, since $\Om$ is closed and $N$ is a
$J$-anti-invariant 2-form with values in $TM$, one easily deduces that $\na
\Om$ is, in fact, a section of the vector bundle $\leftr \La^{0,1}M\otimes
\La^{0,2}M\rightr$, i.e., the following relation holds:
\begin{equation}\label{naJ}
\na_{JX} J = -J(\na_X J), \ \ \forall X \in TM. \end{equation}

\subsection{The Weitzenb{\"o}ck decomposition of the symplectic form}

Let $\delta$ be the codifferential of $g$, acting on (real) $r$-tensors by
$$(\delta T)_{X_1,...,X_{r-1}} = - \sum_{i=1}^{n}(\na_{e_i} T)_{e_i,
X_1,...,X_{r-1}},$$ where $\{e_i\}_{i=1}^{n}$ is any orthonormal frame.
Note that when acting on $r$-forms, $\delta$ is the formal adjoint operator
to the exterior differential $d$; the two operators are then related via
the Hodge star-operator $*_g$: $$ \delta = -*_g \circ *_g {d}.$$ Then, the
operator $\Delta = d\delta + \delta d$, which acts on the smooth sections
of $\La^r M$, is the {\it Hodge-de Rham Laplacian} of $(M,g)$; we also
consider the {\it rough} Laplacian, $\na^*\na$, which for any smooth
$r$-form, $\psi$, gives $\na^*\na \psi : = \delta (\na \psi)$, again an
$r$-form. The two Laplacians are related by a {\it Weitzenb{\"o}ck
decomposition} $$ \Delta = \na^*\na + A(R) ,$$
where $A$ depends linearly on the curvature $R$ of $(M,g)$, see e.g.
\cite{besse}. For example, the Weitzenb{\"o}ck decomposition on 2-forms
reads as:
\begin{eqnarray} \label{wtz2}
\Delta \psi -\na^* \na \psi &=& [{\rm Ric}(\psi \cdot,\cdot) - {\rm
Ric}(\cdot, \psi\cdot)] - 2 {R}(\psi) \\ \nonumber &=&
\frac{2(m-1)}{m(2m-1)}s \psi -2{W}(\psi) \\\nonumber & & +
\frac{(m-2)}{(m-1)}[{\rm Ric}_0(\psi\cdot,\cdot) - {\rm Ric}_0(\cdot,
\psi\cdot)].
\end{eqnarray}

The symplectic form $\Omega $ is a real harmonic 2-form of type $(1,1)$
with respect to any compatible almost-K{\"a}hler metric $(g,J)$, i.e.
$$\Om(J\cdot, J\cdot)= \Om(\cdot,\cdot), \; \; {d} \Om = 0 \mbox{ and }
\delta \Om = 0.$$ Applying relation (\ref{wtz2}) to $\Om$, we obtain
\begin{equation} \label{wtz2Om} \na^* \na \Om = 2{R} (\Om) - [{\rm
Ric}(J\cdot, \cdot) - {\rm Ric}(\cdot, J\cdot )] \; .
\end{equation}
Formula (\ref{wtz2Om}) is a measure of the difference of the two types of
Ricci forms on an almost-K{\"a}hler manifold: \begin{equation} \label{r*r}
\rho^* - \rho = \frac{1}{2}(\na^* \na \Om ). \; \end{equation} Taking the
inner product with $\Om$ of the relation (\ref{r*r}), we obtain the
difference of the two types of scalar curvatures: \begin{equation}
\label{s*s} s^* - s = |\na \Om|^2 = \frac{1}{2} |\na J|^2 \; .
\end{equation} Formulae (\ref{r*r}) and (\ref{s*s}) can be interpreted as
``obstructions'' to the (local) existence of a strictly almost-K{\"a}hler
structure $J$, compatible with a given metric $g$. Indeed, if we denote by
$P$ the curvature type operator acting on $\La^2M$ by $$P(\psi) =
\frac{2(m-1)}{m(2m-1)}s \psi -2{W}(\psi),$$ then, by (\ref{wtz2}) and
(\ref{s*s}), $(P(\Om), \Om ) = -|\na \Om|^2\le 0$. This shows that
Riemannian metrics for which ${P}$ is semi-positive definite do not admit
even locally compatible strictly almost-K{\"a}hler structures
\cite{her-lamoneda}. Note that the latter curvature condition is implied by
the non-negativity of the {\it isotropic} sectional curvatures. The
following compilation of results is an illustration of this criterion of
non-existence: \begin{prop} Let $(M,g)$ be a Riemannian manifold of
dimension $n=2m$, which satisfies one of the following conditions:
\begin{enumerate}
\item[(i)] $M$ is four-dimensional and oriented, and $g$ is anti-self-dual
metric of non-negative scalar curvature; \item[(ii)] $(M,g)$ is a
conformally-flat manifold of non-negative scalar curvature; \item[(iii)]
$(M,g)$ is a symmetric space of compact type. \end{enumerate} Then, $(M,g)$
does not admit even locally defined strictly almost-\ka structures
{\rm(}which are also compatible with the fixed orientation in the case {\rm
(i))}. \end{prop}

\noindent
{\it Remark.} The situation dramatically changes
if one considers symmetric spaces of non-compact type, see Proposition
\ref{prop:adm2} below. $\Diamond$

\section{Cartan-K{\"a}hler Theory: Further obstructions} In this section,
which is inspired from \cite{Arm0} (see also \cite{AA}), we shall consider
the problem of local existence of Einstein almost-\ka metrics. By Darboux's
theorem, it is enough to look for Einstein metrics $$g=
\sum_{i,j=1}^{2m}g_{ij}(x)dx_i\otimes dx_j$$ compatible with the standard
symplectic form
$$\Om_0 = \sum_{i=1}^{m} dx_i\wedge dx_{m+i}$$ on an open set of ${\mathbb
R}^{2m} = \{{x}=(x_1,...,x_{2m}) \}$. By a result of DeTurck and Kazdan
\cite{De-Ka}, any Einstein metric is real analytic (in suitable local
coordinates), so that it is also natural to require the analyticity of the
components $g_{ij}({x})$ of the metric.

The Ricci curvature can be thought as a non-linear second-order
differential operator
$${\rm Ric} : {\bf AK}(\Om_0) \ni g \to {\rm Ric}(g) \in \Gamma (S^2)$$
acting on the space of $\Om_0$-compatible almost-\ka metrics. Here and
henceforth, $S^l$ stands for the bundle of symmetric $l$-tensors on
${\mathbb R}^{2m}$. For a fixed $J\in {\bf AK}(\Om_0)$, let $E= \leftr
S^{0,2} \rightr$ denote the corresponding bundle of $J$-anti-invariant
elements of $S^2$ (see section 2.1). Then the space of smooth sections of
$E$ is naturally identified with the tangent space of ${\bf AK}(\Om_0)$ at
$(g,J)$, and the {\it principal symbol} $\sigma({\rm Ric})$ of the Ricci
operator is a bundle-map
$$\sigma({\rm Ric}) : S^2\otimes E
\longrightarrow S^2$$ given by
$$ \sigma({\rm Ric})(C)_{X,Y} =
\frac{1}{2}\sum_{i=1}^{2m}\Big[C_{e_i,X, e_i, Y} + C_{e_i,Y,e_i,X} -
C_{e_i,e_i,X,Y} - C_{X,Y,e_i,e_i}\Big]$$ for any section $C$ of $S^2\otimes
E$ (see e.g. \cite{besse}). One can easily check that this map is
surjective. This shows that for any symmetric $(2m\times 2m)$-matrix $r_0$
with constant coefficients there are metrics $g$ which are compatible with
$\Om_0$ and such that, at a given point $x_0 \in {\mathbb R}^{2m}$, we have
${\rm Ric}(g)_{x_0} = r_0$; in other words, at any given point $x_0$ and
for any $r \in \Gamma(S^2)$, one can always find {\it point-wise} (or {\it
algebraic}) solutions to the equation ${\rm Ric}(g)_{x_0} = r_{x_0}$.
Similarly, at any point $x_0$ and for any value $\la$, one can find a
metric $g\in {\bf AK}(\Om_0)$ for which ${\rm Ric}(g)_{x_0} = \la g_{x_0}$.
The next question is whether or not a point-wise solution can be extended
to a real analytic solution defined in a neighborhood of the point. A
coordinate-free realization of the Cauchy-Kowalewski theorem, known as the
Cartan-\ka theory, gives a method of answering this question and we refer
the reader to \cite{bryantetal} for details and references relevant to this
method. We recall here that the basic idea behind Cartan-\ka theory is a
simple one: one builds up order by order real analytic solutions to the
given system of differential equations, by using the point-wise solutions.

Although Cartan-\ka theory is usually used as a tool of proving existence
results, in \cite{Arm0} a different aspect of the theory is emphasized. It
also provides a method to find non-obvious conditions --- which we shall
refer to in this paper as ``obstructions'', but which are called
``torsion'' in \cite{bryantetal} --- that solutions of a differential
equation must satisfy. These obstructions encapsulate in an invariant
manner the elementary fact that derivatives in different coordinate
directions must commute. We shall write down in an explicit way such local
obstructions to finding critical almost-\ka metrics and use them to prove
the integrability of the almost-complex structure in the {compact} case,
provided that certain additional curvature conditions are satisfied. Thus,
Cartan-\ka theory will be used to prove {\it non-existence} results. On the
other hand, because of algebraic difficulties in applying the theory to its
end, no general existence result has been proven so far.

Suppose that one wishes to apply the Cartan-\ka theory (as described in
\cite{gasqui}) in order to prove that strictly almost-\ka Einstein metrics
exist. We have already seen that one can always find examples of 2-jets of
compatible metrics which satisfy ${\rm Ric}(g) = \la g$; the next question
is whether or not one can find algebraic examples of 3-jets satisfying this
equation and its first derivative --- i.e. whether or not one can find
algebraic solutions to the first prolongation of the problem. The symbol of
the first prolongation, $\sigma_1({\rm Ric})$, is a bundle map
$$\sigma_1({\rm Ric}) : S^3 \otimes E \to T^*\otimes S^2. $$ By calculating
the dimension of the image of $\sigma_1({\rm Ric})$, one sees that this map
is not onto. To see this directly, recall that the {\it differential
Bianchi identity} implies $\delta {\rm Ric} = -\frac{1}{2}ds,$ where, we
recall, $\delta$ is the co-differential of $g$ and $s ={\rm trace}_g({\rm
Ric})$ is the scalar curvature. If we can extend a 2-jet solution of the
Einstein equation to a 3-jet solution, we must have $d\la =0$; we thus have
found an obstruction to extending the 2-jet solution to a 3-jet solution,
which simply tells us that an Einstein metric has constant scalar
curvature. If we denote by ${\frak b} : T^*\otimes S^2 \to T^*$ the
equivariant map $${\frak b}(C)_X =
\sum_{i=1}^{2m}(C_{e_i,e_i,X} - \frac{1}{2}C_{X,e_i,e_i}),$$ then we have
the exact sequence: $$0\to S^3\otimes E \buildrel {\sigma_1({\rm
Ric})}\over\longrightarrow T^*\otimes S^2 \buildrel {\frak b}
\over\longrightarrow T^* \to 0.$$ This tells us that this is the only such
obstruction.

Exactly the same obstruction arises for the more general problem of finding
Einstein metrics when one does not insist that the metric is compatible
with the symplectic form. It turns out that if one is looking for metrics
not necessarily compatible with a symplectic form, then there are no
further obstructions and the Cartan-\ka theorem allows one to conclude that
there is a germ of local solutions to the general Einstein equations.
\begin{theo}\label{th1} \cite{gasqui} Let $R_0$ be a given algebraic
curvature tensor at the origin in ${\mathbb R}^n, n\ge 3$ {\rm (see section
3.2)}, and let $g_0$ be an algebraic {\rm (}Riemannian{\rm )} metric such
that the Ricci tensor of $R_0$ with respect to $g_0$ is equal to $\la g_0$
for some constant
$\la$. Then there is a real analytic metric $g$ defined in a neighborhood
of the origin, such that
\begin{enumerate}
\item[(i)] $ g_{x=0} = g_0$ and $(R^g)_{x=0} = R_0$; \item[(ii)] ${\rm
Ric}(g) = \la g$.
\end{enumerate}
\end{theo}

Turning back to our initial problem of finding Einstein almost-\ka metrics,
one should notice that if we could extend the above result to the case of
metrics compatible with $\Om_0$, then strictly almost-\ka solutions would
automatically exist by taking any 2-jet solution $(g_0,J_0)$ of the
Einstein equations, for which the curvature $R_0$ does not belong to the
space of k\"ahlerian algebraic tensors with respect to $J_0$; for instance,
just by calculating dimensions, one sees that there are algebraic 2-jet
solutions to the Einstein equation such that
$W''_0\neq 0$ (see section 3.2).

One then naturally wonders, are there any higher obstructions to extending
$\ell$-jet algebraic solutions of the Einstein equations to $(\ell+1)$-jet
solutions, but this time in the setting of $\Om_0$-compatible metrics?
Unfortunately, it turns out that there are. Indeed, letting $\sigma_2({\rm
Ric})$ be the symbol of the second prolongation, it turns out \cite{Arm0}
that the sequence
$$ 0\to S^4\otimes E \buildrel \sigma_2({\rm Ric})\over\longrightarrow
S^2\otimes S^2 \buildrel\sigma_1({\frak b})\over\longrightarrow T^*\otimes
T^* \to 0 $$ is not exact. By calculating the dimension of the image of
$\sigma_2({\rm Ric})$ one can see that there must be an equivariant
bundle-map ${\frak c} : S^2 \otimes S^2 \to {\mathbb R}$ such that $$ 0\to
S^4\otimes E \buildrel \sigma_2({\rm Ric})\over\longrightarrow S^2\otimes
S^2 \buildrel {\sigma_1({\frak b})\oplus {\frak c}}\over\longrightarrow
T^*\otimes T^*\oplus {\mathbb R} \to 0
$$ is exact. Thus, there is some obstruction to extending $3$-jet solutions
of the Einstein equations in the almost-\ka setting.

The above theory does not lead in a particularly simple way to finding the
obstruction in an explicit form. It merely tells us that there is such an
obstruction and that we have to examine the 4-jet of the metric to find it
--- even though eventually the obstruction takes the form of a condition on
the 3-jet.

The explicit calculations for obtaining the obstruction in most general
form have been carried out in \cite{tedi3} (in the four dimensional case)
and have been later extended in \cite{ADM} for the case of higher
dimensional almost-\ka manifolds; as a final result, we obtain a general
identity which holds for any almost-\ka manifold. In particular cases this
identity was found in \cite{AA,Arm0,oguro-sekigawa,sekigawa}. To state our
result, we introduce \begin{equation}\label{phi} \phi(X,Y) = (
\na_{JX}\Om,\na_Y \Om ),
\end{equation}
which, by (\ref{naJ}), is a semi-positive definite $(1,1)$ form. We then
have \begin{prop}\label{o1}\cite{ADM} For
any almost-\ka structure $(g, J, \Om)$ the following relation holds {\rm
(}notations of section 3{\rm )}: \begin{eqnarray} \label{laplace1} \delta
\big(J \delta (J {\rm Ric}'')\big) &=& - \frac{1}{4}\Delta (|\na \Om|^2) +
2 \delta \big(( \rho^* , \na_{\cdot} \; \Om )\big) + \frac{1}{2}|{\rm
Ric}''|^2 - |(\rho^*)''|^2 \\ \nonumber & & - 2|W''|^2 -
\frac{1}{4}|(\na^*\na \Om)'|^2 - \frac{1}{4}|\phi|^2 + ({\rho}, \phi) -
({\rho}, \na^*\na \Om) \; . \end{eqnarray} In the case of an Einstein
manifold, the above relation reduces to
\begin{equation*}
8\delta \big(( \rho^* , \na_{\cdot} \; \Om )\big) - \Delta (|\na \Om|^2) =
8|W''|^2 + 4|(\rho^*)''|^2 + |(\na^*\na \Om)' |^2 + |\phi|^2 +
\frac{s}{n}|\na \Om|^2;
\end{equation*}
this condition on the 3-jet of an Einstein, almost-\ka manifold is
precisely the obstruction referred to earlier. \end{prop}
Suppose $M$ is a compact Einstein, almost-\ka manifold of positive or zero
scalar curvature. After integrating the above relation over the manifold,
we obtain the following result due to Sekigawa.
\begin{theo}\label{thm:sekigawa}\cite{sekigawa} A compact Einstein,
almost-\ka manifold of non-negative scalar curvature is necessarily \ka.
\end{theo}
More generally, (\ref{laplace1}) provides a non-trivial relation on the
3-jet of any critical metric, which can be also viewed as an obstruction
for the solvability of Equation (\ref{eq:1}). In the case of a (connected)
critical almost-\ka 4-manifold whose curvature component $W''$ identically
vanishes, this relation implies a strong maximum principle for the Nijenhuis tensor $N$ \cite{AA}: if $N$ vanishes at one point, then $N$ is identically ze
ro on $M$. It follows that on any (connected) compact critical strictly
almost-\ka 4-manifold with $W''=0$, the complex rank 2 bundle $L=
\La^{0,1}M \otimes \La^{0,2}M$ has a nowhere vanishing (real) smooth
section, the Nijenhuis tensor $N$ (see section 3.3). As observed by J.
Armstrong \cite{Arm2}, the latter in turn implies strong topological
consequences via the Chern-Weil-Wu formulae:
$$0= c_2(L)(M) = (2c_1^2 + c_2)(M) = (5\chi + 6 \sigma)(M),$$ where
$c_2(L), c_1$ and $c_2$ are the Chern classes of the (rank 2) complex
bundles $L$ and $TM$, while $\chi$ and $\sigma$ are the Euler
characteristic and the signature of $M$. Combining with the Hitchin-Thorpe
inequality $(2\chi \pm 3\sigma)(M) \ge 0$ in the Einstein case
\cite{hitchin1}, we obtain the following integrability result.
\begin{theo}\label{th-aa}\cite{AA} Let $(g,J)$ be a critical almost-\ka
metric on a compact symplectic 4-manifold $(M,\Om)$. Suppose that $W''=0$
{\rm(}see section 3.2{\rm)}. Then $(g,J)$ is \ka provided that $(5\chi +
6\sigma)(M) \neq 0$. In particular, a compact Einstein, almost-\ka
4-manifold is \ka if and only if $W''=0$.
\end{theo}

One might further ask if there are any other obstructions to finding
almost-\ka Einstein metrics. For that we consider the third order
differential operator
$$O_1(g) = \Delta (|\na \Om|^2)
-8\delta \big(( \rho^* , \na_{\cdot} \; \Om )\big) + 8|W''|^2 +
4|(\rho^*)''|^2 + |(\na^*\na \Om)' |^2 + |\phi|^2 + \frac{s}{n}|\na \Om|^2$$
and try to solve the modified system
$$ {\rm Ric}(g) - \la g =0, \; \; O_1(g)=0.$$ J. Armstrong \cite{Arm0}
showed that
there are at least two further obstructions arising in this manner. The
first one is an obstruction to lifting the 4-jet to the 5-jet solution of
the above system, which in the case of 4-manifolds corresponds to a
bundle-map
$${\frak d} : S^3M\otimes S^2M \oplus S^2M \to \La^{1,1}M\cong \La^{1,1}_0M
\oplus {\mathbb R}.$$ Finding the above mentioned obstructions in an
explicit form seems to be a difficult task, even in the case of
4-manifolds. However, along the lines of \cite{Arm0}, one sees that a
certain component of ${\frak d}$, $${\frak d}_{\mathbb R} : S^2M\otimes
S^2M \to {\mathbb R} ,$$ should correspond to a relation expressing the
term $(\na^*\na N , N)$ with lower jets of $(g,J)$, where $N$ is the
Nijenhuis tensor of $J$, also identified to $\na \Om$ via the metric, see
section 3.3. We can write down the relation corresponding to ${\frak
d}_{\mathbb R}$ by specifying the Weitzenb{\"o}ck decomposition of
$TM$-valued 2-forms for the particular section $N$. This is done in
\cite{tedi3} (see also \cite{AAD1}).
\begin{prop} \label{o2}\cite{tedi3}
For any almost-\ka 4-manifold $(M,g,J,\Om)$ the following relation holds
\begin{eqnarray*} \label{laplace2}
\delta \big(J \delta (J {\rm Ric}'')\big) &=& - \frac{1}{4}\Delta (|\na
\Om|^2) - \delta \delta ({\rm Ric}'') - \delta \big(( \rho^* , \na_{\cdot}
\; \Om )\big) + \frac{1}{2}|{\rm Ric}''|^2 + |(\rho^*)''|^2 \\ \nonumber &
& - \frac{1}{2}|(\na^{2}\Om)^{\rm sym}_0|^2 - \frac{3}{32}|\na \Om|^4 -
|W''|^2 - \frac{s}{4} |\na \Om|^2 - ({\rho}, \phi) \; ,
\end{eqnarray*}
where $(\na^2 \Om)^{\rm sym}_0$ is the image of $\na^2 \Om$ under the
bundle map $ (T^*M)^{\otimes 2} \otimes \Lambda^2 M \rightarrow S^2_0 M
\otimes \Lambda^2 M$, which acts as the identity on the second factor and
in the first factor takes the trace free part of the symmetrization.
\end{prop}
Combining the two obstructions given by Propositions \ref{o1} and \ref{o2},
and taking into account that in $4$-dimension $(\nabla^* \nabla \Om)' =
\frac{|\na \Om|^2}{2} \Om$ and $|\phi|^2 = \frac{|\na \Om|^4}{4}$, one
derives the following integrability result.

\begin{theo}\label{th-tedi} \cite{tedi3} Let $(M, \Om)$ be a compact,
4-dimensional symplectic manifold with $(c_1 \cdot [\Om])(M) \geq 0$. Let
$(g,J)$ be an $\Om$-compatible, critical almost-\ka metric and assume that
one of the following is satisfied:

(a) the scalar curvature $s$ is non-negative;

(b) the scalar curvature $s$ is constant. \newline Then $(g,J)$ is
necessarily a K\"ahler structure. \end{theo}

\noindent {\it Remarks.} 1. Note that if the scalar curvature is
non-negative, then the condition $(c_1 \cdot [\Om])(M) \geq 0$ is
automatically satisfied because of (\ref{eq:blair}). For the case (b), it
is {\it a priori} open the possibility that the scalar curvature is a
negative constant, although {\it a posteriori}, the result shows that this
cannot be the case, again because of (\ref{eq:blair}). It would be
interesting to see if the result is still valid without any of the
conditions (a) or (b).

2. The above theorem should be seen in the vein of the recent
classification results of \cite{lalonde-mcduff,Liu,ohta-ono}; accordingly,
any compact symplectic 4-manifold which admits a metric of positive scalar
curvature is deformation equivalent to a rational or a ruled \ka surface.
It is also known \cite{hitchin,sung,yau1} that these manifolds admit \ka
structures of positive scalar curvature.

3. As observed in \cite{boyer-galicki}, the integrability results stated in
theorems \ref{thm:sekigawa} and \ref{th-tedi} hold true for compact
almost-\ka orbifolds. $\Diamond$

\section{Examples of critical almost-\ka metrics} In this section we review
some (mostly non-compact) examples of critical (Einstein) strictly
almost-\ka metrics. We first make the remark that given such a manifold we
can take its product with any \ka (Einstein) manifold to produce examples
of higher dimension. We will be therefore interested in irreducible (i.e.
non-product) examples.

\subsection{The B\'erard Bergery construction} Starting from a critical
(Einstein) strictly almost-\ka manifold $(M^{2m},g, J,\Om)$, there is a
method of producing higher dimensional non-product examples, due to L.
B\'erard Bergery \cite{berard-bergery} (see also \cite[Th.9.129]{besse}).
We take an $S^1$-bundle $P$ over $M$, whose curvature is a multiple of
$\Om$. Then the almost-\ka structure on $M$ induces a $K$-contact structure
on $P$ (see e.g. \cite{boyer-galicki}) and further, a critical strictly
almost-\ka structure on the ``symplectisation'' $N^{2(m+1)} =P\times
{\mathbb R}$ of $P$ (see also \cite{Le-Wang}); the metric on $N$ is
complete provided that $M$ is compact. Furthermore, if $M$ is Einstein,
then the strictly almost-\ka metric on $N$ is also Einstein (with scalar
curvature smaller than the one of $M$). The reader is referred to
\cite{berard-bergery} for more details on the construction.

By applying this procedure to the next examples, one can construct many
strictly almost-\ka critical/Einstein manifolds.

\subsection{Twistorial examples}
For a given oriented Riemannian 4-manifold $(M,g)$, the set of all
almost-Hermitian structures compatible with the metric and the orientation
can be naturally identified with the set of sections of the sphere bundle
$S(\La^+ M)$, by identifying each positively oriented almost-Hermitian
structure $(g,J)$ with its (normalized) fundamental 2-form $\frac{1}{\sqrt
2}g(J\cdot,\cdot)$. The total space ${ Z}$ of $S(\La^+ M)$ is called
(positive) {\it twistor space} of $(M,g)$. Thus, $Z$ can be viewed as an
$SO(4)/U(2)$ fibre-bundle associated to the canonical principal
$O(4)$-bundle $P$ of $(M,g)$, i.e., ${Z}\cong
P\times_{{O(4)}}({SO(4)/U(2)})$. Denote by $p: {Z} \longmapsto M$ the
natural projection; the vertical distribution ${\cal V}= $Ker$(p_*)$
inherits a canonical complex structure $J^{\cal V}$ coming from the natural
complex structure of the fibre ${\mathbb CP}^1 ={SO(4)/U(2)}$. Moreover,
the Levi-Civita connection $\nabla$ of $(M,g)$ induces a splitting $T{Z}
={\cal H}\oplus {\cal V} $ of the tangent bundle of ${ Z}$ into horizontal
and vertical components, so that ${\cal H}\cong p_{*}^{-1}TM $ acquires a
tautological complex structure $J^{{\cal H}}$ given by $J^{\cal H}_{x,j} =
j$ for $x\in M$ and $j\in p^{-1}(x)\cong {\Bbb CP}^1$. Following \cite{ES},
we define an almost-complex structure ${J}$ on $Z$, by
\begin{eqnarray}\nonumber
{J} &=& J^{{\cal H}} - J^{{\cal V}}. \end{eqnarray} Using the splitting of
the tangent bundle of ${Z}$ into horizontal and vertical components, one
also determines a family of Riemannian metrics $h$ on ${ Z}$
$$h_t = p^*g + tg_{{\Bbb CP}^1}, \ t>0, $$ where $g_{{\Bbb CP}^1}$ is the
standard metric (of constant curvature $1$) on the fibre ${\Bbb CP}^1$.
Clearly, $h_t$ is compatible with
${J}$. We then have the following
\begin{theo}\cite{DM} \label{prop:dm} The almost-Hermitian 6-manifold
$(Z,h_t,{J})$ is almost-\ka if and only if $(M,g)$ is Einstein,
anti-self-dual 4-manifold of negative scalar curvature $s=-\frac{12}{t}$.
Moreover, in this case $(h_{-\frac{12}{s}},J)$ is strictly almost-\ka
structure whose curvature tensor belongs to ${\cal R}_2(Z)$; in particular,
it is a critical, strictly almost-\ka metric on $Z$. \end{theo}

\vspace{0.2cm} \noindent {\it Example 2.} The only known examples of {\it
compact} Einstein self-dual 4-manifolds of negative scalar curvature are
compact quotients of the real hyperbolic space ${\mathbb RH}^4$, or of the
complex hyperbolic space ${\mathbb CH}^2$, where the quotients of ${\mathbb
CH}^2$ are considered with the non-standard orientation (so that the
induced metric is {\it anti-self-dual}). Let $(M,g)$ be such a quotient and
by rescaling the metric, suppose that the scalar curvature is equal to
$-12$. According to Theorem \ref{prop:dm} the twistor space $(Z,h_1,{J})$
is then a {compact} locally-homogeneous critical almost-\ka 6-manifold
whose curvature belongs to ${\cal R}_2(Z)$. In fact, it turns out that $Z$
is an example of a compact, locally 3-symmetric almost-\ka 6-manifold in
the sense of \cite{gray1}. $\Diamond$

\vspace{0.2cm} \noindent {\it Remarks.} 1. Note that the fundamental group
of $Z$ is equal to that of $M$. When $M$ is a compact quotient of ${\mathbb
RH}^4$, it follows that the smooth 6-manifold underlying the twistor space
${ Z}$ admits no K\"ahler structures at all \cite{carlson-toledo}. The
analogous assertion for the twistor space of a compact quotient of
${\mathbb CH}^2$ is an open and intriguing problem \cite{kotschick-private}.

2. Theorem \ref{prop:dm} was generalized in \cite{AGI} for twistor spaces
of quaternion-\ka manifolds of negative scalar curvature; considering
compact quotients of locally-symmetric quaternion-\ka $4k$-manifolds of
non-compact type, one obtains compact examples of critical strictly
almost-\ka manifolds of dimension $n=4k+2$. $\Diamond$

\subsection{Homogeneous examples via \ka geometry} The examples given in
this section are studied in \cite{ADM}, and rely on the following two
simple lemmas. For their proofs, we refer the reader to \cite{ADM}.
\begin{Lemma}\label{lem:adm1} Let $(M,g,I)$ be a K{\"a}hler manifold whose
Ricci tensor {\rm (}considered as a symmetric endomorphism of the tangent
bundle via the metric $g${\rm )} has constant eigenvalues $\la < \mu$.
Denote by $E_{\la}$ and $E_{\mu}$ the corresponding $I$-invariant
eigenspaces and define an almost-complex structure ${J}$ by
${J}_{|E_{\la}}= I_{|E_{\la}}; \ {J}_{| E_{\mu}}= -I_{| E_{\mu}}$. Then
$(g, {J})$ is an almost-K{\"a}hler structure; it is K{\"a}hler if and only
if $(M,g,I)$ is locally product of two K{\"a}hler--Einstein manifolds of
scalar curvatures $\la$ and $\mu$, respectively. \end{Lemma}

\begin{Lemma}\label{lem:adm2}
With the notations from the previous lemma, let $(g,I)$ be a K{\"a}hler
structure whose Ricci tensor has constant eigenvalues, and consider the
1-parameter family of metrics $g_t=g_{|_{E_{\la}}} + t g_{|_{E_{\mu}}}, \ \
t>0, $ where $g_{|_{E_{\la}}}$ {\rm (}resp. $g_{|_{E_{\mu}}}${\rm )} is the
restriction of $g$ to the eigenspace $E_{\la}$ {\rm (}resp. to
$E_{\mu}${\rm )}. Then, for any $t>0$, the metric $g_t$ is K{\"a}hler with
respect to $I$, almost-K{\"a}hler with respect to ${J}$ and has the same
Ricci tensor as the metric $g = g_1$. \end{Lemma}

\noindent Lemma \ref{lem:adm2} shows that one can normalize any \ka metric
with Ricci tensor having two distinct constant eigenvalues to one whose
Ricci tensor has eigenvalues equal to $-1,0$ or $+1$. In particular, we get

\begin{cor} \label{cor:adm1} On a complex manifold $(M^{2m},I)$ there is a
one-to-one
correspondence between K{\"a}hler metrics $(g,I)$ with Ricci tensor having
constant eigenvalues $\la < \mu$ with $\la\mu >0$ and K{\"a}hler-Einstein
metrics
$({\widetilde g}, I)$ of scalar curvature $2m\la$ carrying an
almost-K{\"a}hler structure $J$ which commutes with and differs from $\pm
I$; in this correspondence ${J}$ is compatible also with $g$ and coincides
{\rm (}up to sign{\rm )} with the almost-K{\"a}hler structure defined in
Lemma \ref{lem:adm1}; moreover, ${J}$ is integrable precisely when $(g,I)$
{\rm (}and $({\widetilde g},I)${\rm )} is locally product of two
K{\"a}hler-Einstein metrics. \end{cor}

\noindent Combining Corollary \ref{cor:adm1} with Theorem
\ref{thm:sekigawa}, it follows that a {\it compact} K{\"a}hler manifold
whose Ricci tensor has two distinct constant {\it positive} eigenvalues is
always locally the product of two K{\"a}hler-Einstein manifolds. In fact,
using Proposition \ref{o1} for the almost-\ka structure $(g,{J})$, one can
refine the latter result:
\begin{theo}\label{thm:adm1} \cite{ADM} Let $(M,g,I)$ be a compact
K{\"a}hler manifold whose Ricci tensor has two distinct constant {\rm
non-negative}
eigenvalues $\la$ and $\mu$.
Then the universal cover of $(M,g,I)$ is the product of two simply
connected K{\"a}hler-Einstein manifolds of scalar curvatures $\la$ and
$\mu$, respectively. \end{theo}
The above theorem can be equally seen as an integrability result for the
almost-\ka metric $(g,{J})$ given by Lemma \ref{lem:adm1}; it is thus
rather disappointing from the point of view of the search for critical
strictly almost-\ka metrics on compact manifolds. In dimension four, the
situation is even more rigid, since we can further improve the above
integrability result, by using the Kodaira classification of compact
complex surfaces and some consequences of the Seiberg-Witten theory,
recently established in \cite{kotschick0,leung,petean}.
\begin{theo}\label{thm:adm2} \cite{ADM} Let $(M,g,I)$ be a compact
K{\"a}hler surface whose Ricci tensor has two distinct constant
eigenvalues. Then one of the following alternatives holds:
\begin{enumerate} \item[{\rm (i)}] $(M,g,I)$ is locally symmetric, i.e.,
locally is the product of Riemann surfaces of distinct constant Gauss
curvatures; \item[{\rm (ii)}] if $(M,g,I)$ is not as described in (i), then
the eigenvalues of the Ricci tensor are both negative and $(M,J)$ must be a
minimal surface of general type with ample canonical bundle and with even
and positive signature. Moreover, in this case, reversing the orientation,
the manifold would admit an Einstein, strictly almost-K{\"a}hler metric.
\end{enumerate}
\end{theo}
At this time, we do not know if the alternative (ii) may really hold. As
pointed out, if it does, it also provides a counter-example to the
four-dimensional Goldberg conjecture.

\vspace{0.2cm} While Lemma \ref{lem:adm1} does not seem to be a
particularly useful tool of constructing {\it compact} examples of critical
almost-\ka metrics, it does provide complete ones, by considering
irreducible homogeneous \ka manifolds having Ricci tensor with two distinct
eigenvalues. Indeed, according to the structure theorem for homogeneous \ka
manifolds \cite{Dor-Nak,V-G-PSh}, there are such examples of any complex
dimension $m\ge 2$. To see this, recall that any homogeneous K{\"a}hler
manifold admits a holomorphic fibering over a homogeneous bounded domain
whose fiber, with the induced K{\"a}hler structure, is isomorphic to a
direct product of a flat homogeneous K{\"a}hler manifold and a simply
connected compact homogeneous K{\"a}hler manifold; in this structure
theorem an important role is played by the Ricci tensor whose kernel
corresponds to the flat factor; thus, when the Ricci tensor is
non-negative, the manifold splits as the product of a flat homogeneous
manifold (corresponding to the kernel of ${\rm Ric}$) and a compact
homogeneous K{\"a}hler manifold (and thus having positive Ricci form).
Considering non-trivial (K\"ahler) homogeneous fibrations over bounded
homogeneous domains, we obtain examples of (deRham) irreducible homogeneous
K{\"a}hler manifolds with two distinct eigenvalues of the Ricci tensor. In
complex dimension 2, there is only one such manifold which is described
below.

\vspace{0.2cm}
\noindent
{\it Example 3.} There is an irreducible homogeneous \ka surface $(M,I,g)$
whose Ricci tensor (after normalization) has eigenvalues
$(-1, 0)$. Following \cite{wall}, $M$ can be written as $$M= \frac{{\mathbb
R}^2\cdot Sl_2({\mathbb R})}{SO(2)},$$ equipped with a left-invariant
K\"ahler structure $(g,I)$. To make the description of the K\"ahler metric
more explicit, one observes that the four-dimensional (real) Lie group
$G={\mathbb R}^2 \cdot Sol_2$ acts on $M$ simply transitively, so that one
can identify $M$ with $G$. The corresponding Lie algebra $\frak{g}$ is then
generated by $X_1,X_2,E_1,E_2$ with
$$[X_1,X_2]=[X_1,E_2]=0; \ [E_1,E_2] = 2E_2; \ [X_2,E_2] = -2X_1; $$
$$[X_1,E_1]=-E_1;\ [X_2,E_1] = X_2. $$
The left-invariant complex structure $I$ is defined by $I(X_1) = X_2; \
I(E_1)=E_2$, while the K\"ahler form is $\Omega_I = X^1\wedge X^2 +
E^1\wedge E^2$. This example is studied in detail in \cite{wall} where it
is referred to as the ${\bf F}_4$-geometry. As a matter of fact, $M$ does
not admit any compact quotients, but it does admit quotients of finite
volume. The corresponding critical strictly almost-\ka structure $(g,J)$
via Lemma 1 is isomorphic to the unique {\it proper 3-symmetric space} in
four dimensions (see \cite{gray1,kow}); it follows by \cite{gray1} that its
curvature tensor belongs to ${\cal R}_2(M)$. The uniqueness of this
example, as well as its coordinate realization, will be discussed in the
next section. $\Diamond$

\vspace{0.2cm} \noindent {\it Remark.} It is erroneously stated in
\cite{ADM} and \cite[Example 3]{ACG} that the homogeneous K\"ahler
surface $M=(SU(2)\cdot Sol_2)/U(1)$ appearing in the
classification of Shima \cite{shima} is deRham irreducible. In
fact, a more careful investigation of the construction in
\cite{shima} shows that any invariant K\"ahler metric $(g,I)$ on
$M=(SU(2)\cdot {Sol}_2)/{U(1)}$ is the K\"ahler product of a
metric of constant Gauss curvature on ${\mathbb CP}^1$ with a
metric of constant Gauss curvature on ${\mathbb CH}^1$. In
particular, the almost-K\"ahler structure $J$ obtained via Lemma~1
is integrable for this example. The issue of existence of local
examples of strictly almost-K\"ahler 4-manifolds with pointwise
constant Lagrangian sectional curvature, raised in \cite{ACG},
should be therefore considered an open problem. $\Diamond$
\vspace{0.2cm}

We are now going to provide (non-compact) examples of deRham irreducible
homogeneous K{\"a}hler manifolds $(M,g,I)$ with Ricci tensor having two
negative eigenvalues $\la< \mu <0$. According to Corollary \ref{cor:adm1},
these will also provide {\it complete} examples of Einstein strictly
almost-K{\"a}hler manifolds. Since ${\rm Ric}$ is negative definite,
$(M,I)$ must be a bounded homogeneous domain; it is a result of Vinberg,
Gindikin and Piatetskii-Shapiro \cite{V-G-PSh} that any such domain has a
realization as a {\it Siegel domain of type II}, i.e, a domain $D = \{
(z,w) \in {\Bbb C}^p\times {\Bbb C}^q : {\rm Im} z - H(w,w) \in
\mathfrak{C} \}$, where
$\mathfrak{C}$ is an open convex cone (containing no lines) in ${\Bbb R}^p$
and $H:{\Bbb C}^q\times {\Bbb C}^q \mapsto {\Bbb C}^p$ is a Hermitian map
which is $\mathfrak{C}$-positive in the sense that $$H(w,w) \in {\overline
{\mathfrak{C}}} -\{0\} \ \ \ \forall w \neq 0.$$ In the particular case
when $q=0$, we have $D= {\Bbb R}^p + i \mathfrak{C}$; such domains are
called {\it Siegel domains of type I} or {\it tube domains}. The following
observation made in \cite{ADM} provides the needed examples:
\begin{prop}\label{prop:adm2} Every irreducible Hermitian symmetric space
of non-compact type which admits a realization as a tube domain carries a
strictly almost-K{\"a}hler structure commuting with the standard K{\"a}hler
structure.
\end{prop}

\noindent
{\it Example 4.} Explicit examples of Hermitian symmetric spaces that admit
strictly almost-K{\"a}hler structures are $$M^{2m} =
\frac{SO(2,m)}{SO(2)\times {SO}(m)}, \ m\ge 3.$$ Indeed, it is well known
that these spaces admit realizations as tube domains, cf. e.g.
\cite{satake}. $\Diamond$

\vspace{0.2cm} \noindent {\it Remarks.} 1. By using Corollary
\ref{cor:adm1}, one can also find {\it non-symmetric} homogeneous
K\"ahler-Einstein manifolds of complex dimension greater than three, which
admit strictly almost-\ka structures. In complex dimension two and three,
however, any bounded homogeneous domain is symmetric \cite{cartan}.
Therefore, according to Corollary \ref{cor:symmetric-spaces} in section 6,
no four-dimensional examples of Einstein, strictly almost-\ka metrics arise
from the above construction.

2. The reader may ask whether {\it compact} Einstein strictly
almost-\ka examples can be found as (smooth) quotients of the
homogeneous examples presented. The answer is negative. Indeed, it
is well-known that if a homogeneous space $M = G/H$ admits a
compact smooth quotient $M/\Gamma$, then $G$ must be unimodular
(see e.g. \cite{ragunathan}); it then follows by
\cite{borel2,hano} that a bounded homogeneous domain admits
compact quotients if and only if it is symmetric. On the other
hand, the maximal connected subgroup $G'$ in $G$ of isometries of
the almost-\ka structure arising from Proposition \ref{prop:adm2}
is a proper closed subgroup of the full group of isometries $G$ of
the Hermitian symmetric space $(M,g)$. To see this, observe that
according to Lemma \ref{lem:adm2} there is at least a
one-parameter family of $G'$-invariant \ka metrics on $M$, while
for an irreducible Hermitian symmetric space the Bergman metric
$g$ is the unique $G$-invariant metric. It is known \cite{borel1}
that no co-compact lattice $\Gamma$ of $G$ can be contained in a
connected proper closed subgroup of $G$, so that the almost-\ka
structure $J$ does not descend to any compact quotient.

3. The first homogeneous examples of Einstein, strictly almost-\ka metrics
have been found by D. Alekseevsky \cite{Ale} (see also
\cite[14.100]{besse}) on certain solvable Lie groups of real dimension $n=
4k, k\ge 6$. The Alekseevsky examples appear in the framework of
quaternion-\ka geometry and differ from the ones arising from homogeneous
\ka manifolds via Corollary \ref{cor:adm1}. Note that the Alekseevsky
spaces do not admit compact quotients either \cite{alekseevsky-cortes}.
$\Diamond$

\subsection{Four-dimensional examples}
The following construction is discussed in \cite{AAD,AAD1} and generalizes
the examples described in \cite{NuP} and \cite{Arm2}.
\begin{prop}\label{prop:AK3}
Let $(\Sigma, g_{\Sigma}, \om_{\Sigma})$ be an oriented Riemann surface
with metric $g_{\Sigma}$ and volume form $\omega_{\Sigma}$, and let $h = w
+ iv$ be a non-constant holomorphic function on $\Sigma$, whose real part
$w$ is everywhere positive. On the product of $\Sigma$ with ${\mathbb R}^2
= \{ (z,t) \}$ consider the symplectic form \begin{equation}\label{om0}
\Om = \om_{\Sigma} - dz\wedge dt
\end{equation}
and the compatible Riemannian metric
\begin{equation}\label{g0}
g = g_{\Sigma} + w dz^{\otimes 2} + \frac{1}{w} (dt + vdz)^{\otimes 2};
\end{equation}
Then, $(g,\Om)$ defines a strictly almost-\ka structure whose curvature
belongs to ${\cal R}_3(M)$; in particular, $g$ is a critical almost-\ka
metric.
\end{prop}
Clearly, $\frac{\partial}{\partial t}$ and $\frac{\partial }{\partial z}$
are two commuting hamiltonian Killing vector fields for the almost-\ka
structure of Proposition \ref{prop:AK3}. Thus, the metric (\ref{g0}) is
endowed with an ${\mathbb R}^2$-isometric action which is {\it
surface-orthogonal}. This suggests that the construction is local in nature
and will not lead to compact examples (see Theorem \ref{theo:AK3glob}
below).

The almost-\ka structure (\ref{om0}-\ref{g0}) is sufficiently explicit to
make it straightforward, though tedious, to check the following assertions:
\begin{enumerate}
\item[$\bullet$] $g$ is Einstein if and only if $(\Sigma, g_{\Sigma})$ is a
flat surface, in which case $g$ is a Gibbons-Hawking self-dual {\it
Ricci-flat} metric
\cite{gibbons-hawking} with respect to a translation-invariant harmonic
function $w$ (see \cite{Arm1,NuP});
\item[$\bullet$] the curvature of $g$ belongs to ${\cal R}_2(M)$
(equivalently, $W'' =0$) if and only if $(\Sigma,g_{\Sigma})$is hyperbolic,
and, in the half-plane realization $({\mathbb RH}^2, \frac{dx^{2} +
dy^{2}}{x^2})$, the holomorphic function $h$ is given by $h= x+ iy$; in
this case (\ref{om0}-\ref{g0}) is homogeneous and isometric to the
four-dimensional proper 3-symmetric space described in Example 3, see
\cite{AAD}. \end{enumerate}

We also note here that an important feature of the above construction comes
from the following observation: at any point where $dw \neq 0$, the \ka
nullity ${D} = \{ TM \ni X : \na_X J = 0 \}$ of $J$ is a two dimensional
subspace of $TM$, which is tangent to the surface $\Sigma$. If we define a
new almost-complex structure $I$ by $$J_{D} = I_{D} , \ \ J_{{D}^{\perp}} =
- I_{{D}^{\perp}},$$ where $D^{\perp}$ is the $g$-orthogonal complement of
$D$, then one immediately sees that $(g,I)$ is a \ka structure. From this
point of view, Proposition \ref{prop:AK3} provides a construction similar
to the one coming from Lemma \ref{lem:adm1}. In fact, at each point, the
Ricci tensor of (\ref{g0}) has eigenvalues equal to $(0,\frac{s}{2})$ and
these are constant precisely when $s$ is constant. If we introduce local
coordinates on $(\Sigma, g_{\Sigma})$, such that $g_{\sigma} =
e^{u}w(dx^{\otimes 2} + dy^{\otimes 2})$ for some function $u(x,y)$, then
the latter condition reads $$ u_{xx} + u_{yy} = (sw) e^{u}$$
for a constant $s$. Generic (local) solutions, $u$, to the above equation
provide {\it non-homogeneous} metrics. Hence there are non-homogeneous
examples satisfying the conditions of Lemma \ref{lem:adm1}.

\vspace{0.2cm}
The next proposition is proved in \cite{ACG}. \begin{prop}\label{prop:wsks}
Let $g_\Sigma$ be a metric on a $2$-manifold $\Sigma$ with volume form
$\omega_\Sigma$, $\beta$ be a $1$-form on $\Sigma$ with
$d\beta=w\,\omega_{\Sigma}$ for an arbitrary holomorphic function $h= w +
iv$ on $\Sigma$, whose real
part $w$ is everywhere positive. Then
\begin{align}\label{ak-metric}
g &= \frac{w}{z}(z^2 g_{\Sigma} + dz^{\otimes 2}) + \frac{z}{w}\Bigl( dt +
\frac{v}{z} dz+\beta\Bigr)^{\otimes 2},\\ \label{ak-form} \Omega &= z
w\omega_{\Sigma} + dz\wedge\Bigl( dt + \frac{v}{z} dz+\beta\Bigr),
\label{ak-omega}\end{align}
defines a critical almost-K{\"a}hler metric which is \ka if and only if the
function $h$ is constant. \end{prop} Note that the almost-\ka metric
(\ref{ak-metric}-\ref{ak-form}) is endowed with a hamiltonian Killing
vector field $K = \frac{\partial}{\partial t}$. Furthermore, one can
directly compute the curvature of (\ref{ak-metric}-\ref{ak-form}). It turns
out that it does not belong to ${\cal R}_3(M)$, showing that the solutions
are different from the ones given by Proposition \ref{prop:AK3}. We thus
can get new examples of self-dual {\it Ricci-flat} strictly almost-\ka
$4$-manifolds. Indeed, let $(g,\Omega)$ be given by
(\ref{ak-metric}-\ref{ak-form}) and suppose that $g_{{\mathbb CP}^1}$ is
the standard metric on an open subset $\Sigma$ of ${\mathbb CP}^1$; let $h=
w + iv$ be a non-constant holomorphic function on $\Sigma$ with positive
real part. If we write $z=r$, we see that the metric $$g = \frac{w}{r} (
dr^{\otimes 2} + r^2 g_{{\mathbb CP}^1}) + \frac{r}{w} (dt+\alpha)^{\otimes
2}$$ is given by applying the Gibbons-Hawking Ansatz \cite{gibbons-hawking}
applied with respect to the harmonic function $w/r$.
This class of Gibbons-Hawking metrics has been studied in~\cite{CaPe:sdc}
and~\cite{CaTo:emh}. The reader is referred to these references for more
information. The observation that these metrics are almost-K{\"a}hler was
made in \cite{ACG}.

\section{Local classification results in dimension 4 and further
integrability results.}

We start this section with the following result recently proven in \cite{AAD1}.
\begin{theo}\label{theo:AK3}
For any connected strictly almost-\ka 4-manifold $(M,g,\Om)$ whose
curvature belongs to ${\cal R}_3(M)$ there exists an open dense subset $U$
with the following property: in a neighborhood of any point of $U$,
$(g,\Om)$ is homothetic to an almost-\ka structure given by Proposition
\ref{prop:AK3}.
\end{theo}

This theorem generalizes some previous results:

\begin{cor} \label{cor:ASD} \cite{Arm2} There are no Einstein
anti-self-dual strictly almost-\ka 4-manifolds. \end{cor} Since any
anti-self-dual \ka surface is necessarily scalar-flat (see e.g.
\cite{gauduchon0}), it also follows

\begin{cor} \label{cor:symmetric-spaces} \cite{Arm2} Let $(M,g)$ be the
four dimensional real hyperbolic space, or the two dimensional complex
hyperbolic space but endowed with the non-standard orientation. Then
$(M,g)$ does not admit even a locally defined compatible almost \ka
structure.
\end{cor}

\noindent {\it Remarks.} 1. To the best of our knowledge, it is still an
open problem whether the complex hyperbolic space ${\mathbb CH}^2$ does
admit a (locally defined) strictly almost-\ka structure $J$ which is
compatible with the canonical metric and with the {\it standard}
orientation. If it does, this would be the first example of a strictly
almost-K\"ahler Einstein, self-dual metric with negative scalar curvature,
and would constitute a major step in the development of the local theory.
In the same vein, it is unknown if ${\mathbb CH}^1 \times {\mathbb CH}^1$
(locally) admits strictly almost-\ka structures compatible with the product
Einstein metric. By contrast, for the other symmetric four-dimensional
spaces the answer is known: applying the curvature criterion for
non-existence of strictly almost-\ka structures, established in section
3.3, one concludes that a strictly almost-\ka structure can possibly exist
only on the space ${\mathbb R}\times {\mathbb RH}^3$. Conversely, Oguro and
Sekigawa showed \cite{oguro-sekigawa2} that this symmetric space does admit
a strictly almost-\ka structure.

2. Any real hyperbolic space $({\mathbb RH}^{2m}, {\rm can}), m\ge 2$ does
not admit even locally defined almost-\ka structures. This result was
proved for $m>3$ in \cite{olszak}. The 6-dimensional version is much more
subtle and was proved only recently in \cite{Arm2}. $\Diamond$

\vspace{0.2cm}
As another immediate consequence from Theorem \ref{theo:AK3} and the
results in \cite{gray1}, we derive

\begin{cor} \cite{AAD} Let $(M,g,J,\Om)$ be a complete, simply connected
strictly almost-\ka 4-manifold whose curvature tensor belongs to ${\cal
R}_2(M)$. Then $(M,g,J,\Om)$ is isometric to the proper 3-symmetric space
described in Example 3. \end{cor}

Since the homogeneous space of Example 3 does not admit any compact
quotient (see \cite{wall}), we also get the following integrability result
originally proven in \cite{ADK}.

\begin{cor} Any compact almost-\ka
4-manifold whose curvature tensor belongs to ${\cal R}_2(M)$ is \ka. \end{cor}

We can ask more generally whether there are {\it compact} strictly
almost-\ka 4-manifolds whose curvature belongs to ${\cal R}_3(M)$. Using
the local structure established in Theorem \ref{theo:AK3} and some further
global arguments, it was also shown in \cite{AAD1} that the answer is
negative.

\begin{theo} \cite{AAD1} \label{theo:AK3glob} Any compact almost-\ka
4-manifold whose curvature tensor belongs to ${\cal R}_3(M)$ is \ka.
\end{theo}

The above integrability result also implies \begin{cor}
There are no critical, strictly almost-\ka structures compatible with a
compact locally-homogeneous Riemannian 4-manifold. \end{cor}

\noindent {\it Proof.} Suppose for contradiction that $(M,g)$ is a compact
locally homogeneous 4-manifold which admits a compatible strictly
almost-\ka structure $J$ such that the Ricci tensor ${\rm Ric}$ is
$J$-invariant.

If $g$ is Einstein, then by a well known result of Jensen \cite{jensen}
$(M,g)$ must be a compact locally symmetric 4-manifold (see also
\cite[Prop.9]{derdzinski} for a different proof). Using a case by case
verification, it is shown in \cite{oguro-sekigawa1} that $(M,g)$ does not
admit strictly almost-\ka structures. We give an alternative proof of this
fact as follows. By Theorem \ref{thm:sekigawa}, we can assume that $(M,g)$
is a compact quotient of a symmetric Einstein 4-manifold of non-compact
type; the latter are ${\mathbb RH}^4$, ${\mathbb CH}^2$ and ${\mathbb
CH}^1\times {\mathbb CH}^1$. By Corollary \ref{cor:symmetric-spaces}, we
can further assume that $(M,g)$ is locally isometric to ${\mathbb CH}^2$ or
${\mathbb CH}^1\times {\mathbb CH}^1$, and that $J$ agrees with the
standard orientation on these spaces. Then, being locally Hermitian
symmetric, $(M,g)$ admits both a \ka structure and a strictly almost-\ka
structure compatible with the same orientation. This is impossible
according to \cite[Th.1]{AD}.

Consider now the case when $(M,g)$ is locally homogeneous but not Einstein.
Since ${\rm Ric}$ is $J$-invariant, it has two distinct (constant)
eigenvalues $\la$ and $\mu$, each one of multiplicity $2$, such that the
corresponding eigenspaces are invariant under the action of $J$; it follows
that $J$ is uniquely determined (up to sign) on the oriented 4-manifold
$(M,g)$, by acting as a rotation of angle $+\frac{\pi}{2}$ on each
eigenspace of ${\rm Ric}$. This shows that $J$ must be invariant under the
transitive group action, i.e., that $(g,J)$ is a locally homogeneous
almost-\ka 4-manifold. In particular, the 2-form $(\rho^*)''$ defined in
section 3.2 is invariant, and therefore of constant length. It is easily
seen that $(\rho^*)''=0$ if and only if the curvature of $(g,J)$ belongs to
${\cal R}_3(M)$ (see e.g. \cite[Lemma 1]{AAD1}). Since we assumed $J$ is
non-integrable, by Theorem {\ref{theo:AK3glob} we conclude that
$(\rho^*)''$ is a nowhere vanishing smooth section of the bundle $\leftr
\La^{0,2}M \rightr$. In four dimensions, the induced action of $J$ on
$\leftr \La^{0,2}M \rightr$ gives an isomorphism between this bundle and
the canonical bundle of $(g,J)$. This shows that $c_1^{\mathbb R} (J)=0$.
Now, by Theorem \ref{th-tedi}, $J$ must be integrable, a contradiction.
$\Box$

\vspace{0.2cm}
We end this section with the following characterization of the metrics
given by (\ref{ak-metric}-\ref{ak-form}).

\begin{theo} \cite{ACG} Let $(M,g,J,\Om)$ be an almost-K{\"a}hler
$4$-manifold with $J$-invariant Ricci tensor and a non-vanishing
hamiltonian Killing vector field $K$. Suppose that the pair $(\bar
g=\mu^{-2}g,I)$ is K{\"a}hler, where $\mu$ is a momentum map for a nonzero
multiple of $K$, and $I$ is equal to $J$ on ${\rm span} (K,JK)$, but to
$-J$ on the orthogonal complement of ${\rm span} (K,JK)$. Then either $J$
is integrable, or $(g,\Om)$ is locally given by Proposition
\ref{prop:wsks}. \end{theo}

\section{Curvature conditions with respect to the Hermitian connection}

For an almost-\ka manifold $(M^{2m},g,J,\Om)$, besides the Levi-Civita
connection $\na$ of $g$, one can also consider the so called {\it
Hermitian} or {\it first canonical} connection (see e.g.
\cite{gauduchon1}), defined by : $${\widetilde{\na}}_X Y = \na_X Y -
\frac{1}{2}J(\na_X J)(Y) . $$ This is the unique connection satisfying the
following three conditions:
$$ {\widetilde{\na}} g = 0, \ \ \ {\widetilde{\na}} J = 0, \ \ \
(T^{{\widetilde{\na}}})^{1,1} = 0 ,$$ where $(T^{{\widetilde{\na}}})^{1,1}$
denotes the $J$-invariant part of the torsion $T^{{\widetilde{\na}}}$ of
${\widetilde{\na}}$, viewed as a real 2-form with values in $TM$. For any
connection that has the first two of the above properties, the
$J$-anti-invariant part of its torsion is canonically identified with the
Nijenhuis tensor $N$ of the almost-complex structure; it follows that among
connections making both $g$ and $J$ parallel, the first canonical
connection has torsion with minimal possible norm, cf. \cite{gauduchon1}.
By Chern-Weil theory, the Chern classes of the manifold are directly
related to the curvature of the connection ${\widetilde {\na}}$. If we
denote by $\widetilde{R}$ the curvature tensor of ${\widetilde {\na}}$,
then $$ \widetilde{\rho}(X,Y) = \sum_{i=1}^m(\widetilde{R}_{X,Y} e_i,
Je_i)$$ is a closed 2-form which is a deRham representative of $2\pi
c_1(M,J)$ in $H^2(M,{\mathbb R})$. Moreover, $\widetilde{\rho}$ is related
to the Ricci forms of $\na $ introduced at the end of section 3.2 by
$$ \widetilde{\rho} = \rho^{*} - \frac{1}{2} \phi = \rho + \frac{1}{2}
(\na^*\na \Om - \phi),$$ where for the second equality we used (\ref{r*r})
and $\phi$ is the $(1,1)$-form given by (\ref{phi}). The {\it Hermitian
scalar curvature} is defined to be
$\widetilde{s} = 2(\widetilde{R}(\Om), \Om)$; by (\ref{s*s}), ${\widetilde
s}$ is related to the scalar and $*$-scalar curvatures by $$ \widetilde{s}
= \frac{1}{2}(s^* + s) = s + \frac{1}{2} |\na \Om|^2.$$ Further, we have
$$ \int_M \widetilde{s} \; dv = \frac{4\pi}{(m-1)!} (c_1 \cdot
[\Omega]^{\wedge (m-1)})(M) = c \; vol(M),$$ where $c$ is the constant $
\frac{4\pi}{(m-1)!}(c_1 \cdot [\Omega]^{\wedge (m-1)})(M)/ vol(M)$. Thus,
almost-\ka metrics which satisfy the conditions \begin{equation}
\label{H-Einstein}
\widetilde{\rho} = \frac{\widetilde{s}}{2m} \Om , \ \ \ \ \widetilde{s} = c
\end{equation}
are natural candidates
for generalizing K\"ahler-Einstein metrics. They could exist only on
compact symplectic manifolds for which the first Chern class $c_1$ is a
multiple of the cohomology class of the symplectic form $[\Om]$.

In the compact case, almost-\ka metrics satisfying (\ref{H-Einstein}) can
be equivalently characterized by the following three conditions:
\begin{equation} \label{H-Einstein2}
{\rm (a)} \ \ \; \widetilde{s} = c, \ \ \ {\rm (b)} \ \ \;
\widetilde{\rho}
\mbox{ is of type} \ (1,1), \ \ \ {\rm (c)} \ \ \; c_1 = \frac{c}{4m\pi}
[\Om] .
\end{equation} It is natural to drop the cohomological condition (c) and
look for compact almost-\ka manifolds satisfying just (a) and/or (b). For
instance, in dimension four, almost-\ka metrics which saturate the new
curvature estimates of LeBrun \cite{lebrun3} must satisfy (a) and (b).

The study of almost-\ka metrics of constant Hermitian scalar curvature is
motivated by some results and questions of Donaldson \cite{donaldson1}. He
uses a moment map approach to show that problems posed by Calabi about
finding canonical \ka metrics on complex manifolds have natural extensions
in the symplectic context.

The first observation is that the Fr\'echet space ${\bf AK}(M, \Om)$ of
$\Om$-compatible metrics admits a formal \ka structure which is preserved
by the action of the connected component of the identity, $Symp_0(M,\Om)$,
of the symplectomorphism group of $(M,\Om)$. Assuming $H^1(M, {\mathbb R})
= 0$, the Lie algebra of this group is identified with the space of smooth
functions on $M$ of zero integral, $C^{\infty}_0(M)$, endowed with the
Poisson bracket with respect to $\Om$. Then, it is shown in
\cite{donaldson1} that
$$\mu : {\bf AK}(M,
\Om) \rightarrow C^{\infty}_0(M), \ \ \ \mu (g) = \widetilde{s}_g - c ,$$
is a moment map for the action of ${Symp}_0(M,\Om)$ on ${\bf AK}(M, \Om)$.
The elements in the zero-set $\mu^{-1}(0)$ of the moment map are
$\Om$-compatible almost-\ka metrics with constant Hermitian scalar
curvature. Having in mind the identification of quotients which holds in
the finite dimensional theory of moment maps (see e.g. \cite{kirwan}), but
which is largely conjectural for infinite dimensional problems, one can
speculate that each ``stable'' leaf of the action on ${\bf AK}(M, \Om)$
induced by the complexified Lie algebra of $Symp_0(M,\Om)$, contains
precisely one $Symp_0(M,\Om)$ orbit of almost-\ka metrics of constant
Hermitian scalar curvature. There is supporting evidence for the validity
of this scenario. Restricted to the space of $\Om$-compatible, {\it
integrable} almost-\ka structures and after an application of Moser's
lemma, Donaldson's approach transforms exactly into the classical set up of
Calabi regarding existence of \ka metrics of constant scalar curvature (in
particular, K\"ahler-Einstein metrics) on compact complex manifolds. The
leaves of the complexified action correspond in this case to deformations
of the \ka form in a fixed cohomology class by \ka potentials. A number of
results in \ka geometry partially confirm the quotient identification in
the integrable case.

Of course, the ultimate goal would be to prove a general existence result
in the symplectic context, but for the beginning it would be interesting to
construct some particular compact examples of strictly almost-\ka manifolds
satisfying some of the conditions in (\ref{H-Einstein2}). At this point,
our knowledge of such examples is limited. Of course, all locally
homogeneous almost-\ka structures have constant Hermitian scalar curvature.
Some satisfy the other conditions of (\ref{H-Einstein2}). For instance, it
is not difficult to check that the Kodaira-Thurston manifold presented in
Example 1 satisfies $\widetilde{\rho} = 0$. Using the curvature computation
of \cite{DM}, the same can be seen to be true for the Davidov-Mu\u{s}karov
twistorial examples presented in section 5.1, cf.~\cite{DM1}.

Finally, it is worth mentioning another research problem suggested by
Donaldson's approach --- the study of the critical points of the squared
norm of the moment map, that is, the $L^2$-norm of the Hermitian scalar
curvature, seen as a functional on the space ${\bf AK}(M, \Om)$. It follows
from the moment map set-up that the critical points of this functional are
precisely the almost-\ka metrics for which the vector field dual to $J d
\widetilde{s}$ is Killing with respect to $g$. Thus, we get a natural
extension of the Calabi extremal \ka metrics \cite{calabi} to the
almost-\ka case. It is interesting to see what parts of the theory of
extremal \ka metrics extend to the symplectic context. For example,
similarly to the \ka case \cite{levin}, one immediately finds an obvious
necessary condition for a compact symplectic manifold $(M,\Om)$ to admit a
compatible {\it extremal} almost-\ka metric $g$ of non-constant Hermitian
scalar curvature: since the connected component of the identity of the
isometry group of $g$ is a compact subgroup of
$Symp_0(M,\Om)$
\cite{lichne},
it follows that $Symp_0(M,\Om)$ contains an $S^1$ which acts with fixed
points. In dimension four, combining this observation with the results in
\cite{mcduff} and \cite{karshon}, one sees that any compact symplectic
4-manifold $(M,\Om)$ possibly admitting an extremal almost-\ka metric of
non-constant Hermitian scalar curvature must be symplectomorphic to a
rational or ruled \ka surface endowed with an $S^1$-isometric (and
therefore holomorphic) action. These are complex surfaces possibly
supporting extremal \ka metrics as well (see e.g.
\cite{calabi,lebrun-simanca,christina1,christina2}).

\vspace{0.2cm} \noindent {\bf Acknowledgement.} We owe thanks to J.
Armstrong, D. Calderbank, P. Gauduchon, D. Kotschick and A. Moroianu for
collaboration on a number of results mentioned in this survey, and to D.
Alekseevsky whose comments were at the origin of the examples presented in
section 5.3. The first author would also like to take this opportunity to
thank the organizers for a wonderful conference.

\end{document}